\documentclass[a4paper,12pt,titlepage]{article}
\usepackage[frenchb]{babel}
\usepackage[latin1]{inputenc}
\usepackage{epsfig}
\usepackage{amsmath,amsfonts,amssymb}
\usepackage{graphicx}
\usepackage{epic}
\usepackage{eepic}
\usepackage{ecltree}
\usepackage{stmaryrd}
\newtheorem{theo}{Theorem}
\newtheorem{lem}{Lemma}

\begin{document}

\begin{center}
{\LARGE Le Cam spacings theorem in dimension two}

\vspace{0.5cm}

Lionel Cucala

{\footnotesize LSP, Université Paul Sabatier and GREMAQ, Université
  Sciences-Sociales, 

Toulouse, France}

\vspace{0.5cm}

{\footnotesize E-mail address: cucala@cict.fr}

\vspace{0.5cm}


\end{center}


\noindent MSC $2000$: 60F05, 62G30. 

\section*{Abstract}

The definition of spacings associated to a sequence of random variables is
extended to the case of random vectors in $[0,1]^2$. Beirlant \& \textit{al}.
(1991) give an alternative proof of the Le Cam (1958) theorem concerning
asymptotic normality of additive functions of uniform spacings in $[0,1]$. I
adapt their technique to the two-dimensional case, leading the way to new
directions in the domain of Complete Spatial Randomness (CSR) testing.

\vspace{0.5cm}

\section{Introduction}

Testing the uniformity of real random variables (r.v.) can be done in several
ways: using Chi-square tests, tests based on the empirical distribution
function (e.d.f.), tests based on spacings ... The latter ones have been
extensively studied (Pyke, 1965) and recommended, for example, for the
analysis of the local renewal structure of a point process (Deheuvels, 1983a). 

In higher dimensions, when dealing with a spatial point pattern $U \in S
\subset \mathbb{R}^d$, one first wishes to know whether it satisfies the CSR
hypothesis: is the spatial process governing $U$ a homogeneous Poisson
process? This question is equivalent to the following: given the number of
points in the pattern (also called events), are these points uniformly and
independently distributed in $S$ (Moller \& Waagepetersen, 2004)? 

I concentrate here on point patterns distributed in rectangles in
$\mathbb{R}^2$, which is similar, after linear transformation of the
coordinates, to testing the uniformity in $[0,1]^2$.

\noindent Most of two-dimensional uniformity tests are either Chi-square tests
or distance-based methods (Cressie, 1993). The first ones depend on the number
and location of the quadrats (cells in which events are counted), whereas the
last ones require numerous simulations. More recently, there has been some
interest in e.d.f.-based methods and extensions of
the Cramer-Von Mises test(Zimmerman, 1993) and the Kolmogorov-Smirnov test
(Justel \& \textit{al}, 1997) to the $[0,1]^2$ case have been established.

\noindent On the other hand, spacings theory, so useful for testing uniformity
on $\mathbb{R}$, remains almost unworked in higher dimensions even if one
may think, as Zimmerman (1993) does, that distances from events to their
nearest neighbours can be viewed as two-dimensional analogues of spacings. We
shall just mention the results of Deheuvels (1983b) and Janson (1987)
concerning the asymptotic distribution of the maximal multidimensional
spacing, i.e. the volume of the largest square (or ball) contained in
$[0,1]^2$ and avoiding every point of the pattern.

\noindent A first application of spacings theory to CSR testing would be to
test both x- and y-coordinates' uniformity using a spacings-based method. The
rejection of either leads to refuse the two-dimensional uniformity hypothesis.
But we can never accept it as a bivariate distribution with uniform marginals
need not be uniform. This makes necessary to take into account the joint
distribution of the x- and y-coordinates.

In this paper, following this idea, I introduce a new notion of
two-dimensional spacings which is related to spacings based on x- and
y-coordinates. This relationship then allows me to derive the limiting
distribution, under the uniformity hypothesis, of a wide family of statistics
based on these spacings. This is done by a direct decomposition method similar
to the one Beirlant \& \textit{al} (1991) used for one-dimensional spacings.  

An application of this asymptotic result is developed by Cucala \&
Thomas-Agnan (2005). Two of these statistics, the variance and the absolute
mean deviation of the two-dimensional spacings, are selected and used in
practice to test for CSR. A multiple procedure is adopted to generalize the
tests to point processes in any domain (not necessarily rectangular). Then the
power of these spacings-based tests is compared to the power of existing tests
using real and simulated data sets: they appear to be inferior for detecting
regularity or clustering but more powerful for detecting certain types of
heterogeneity.

\vspace{-0.25cm}

\section{Spacings in $[0,1]^2$}

\subsection{Definition}

Let $\mathsf{U}=\Big( (U^x_1,U^y_1), \cdots, (U^x_{n-1},U^y_{n-1}) \Big)$ be a
point pattern in $[0,1]^2$.

\noindent Let $\mathsf{U^x}=(U^x_1, \cdots, U^x_{n-1})$ and $\mathsf{U^y}=(U^y_1, \cdots, U^y_{n-1})$.

\vspace{0.2cm}

\noindent $U^x_{(1)} \leq \cdots \leq U^x_{(n-1)}$ are the order statistics
corresponding to $\mathsf{U^x}$.

\noindent $U^y_{(1)} \leq \cdots \leq U^y_{(n-1)}$ are the order statistics
corresponding to $\mathsf{U^y}$.

\vspace{0.2cm}

\noindent Set $U^x_0 = U^y_0 = U^x_{(0)} = U^y_{(0)} = 0$ and $U^x_n =
U^y_n = U^x_{(n)} = U^y_{(n)} = 1$.

\vspace{0.2cm}

\noindent One may define the spacings related to the pattern
\begin{eqnarray*}
\textrm{the x-spacings } & D^{x}_i = U^x_{(i)} - U^x_{(i-1)} , \quad i=1,
\cdots, n \, \, \, ,\\
\textrm{the y-spacings } & D^{y}_j = U^y_{(j)} - U^y_{(j-1)} , \quad j=1,
\cdots, n \, \, \, .
\end{eqnarray*}

\begin{figure}
\begin{center}
\includegraphics[scale=0.5]{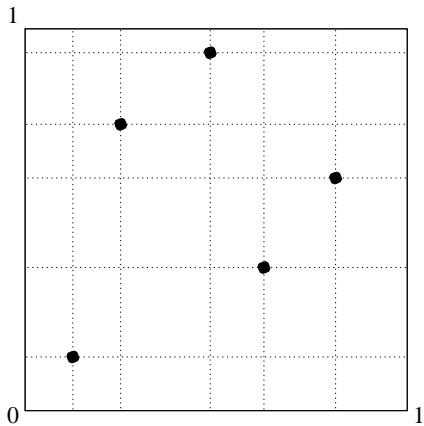}
\caption{Two-dimensional spacings} \label{figzone}
\end{center}
\end{figure}

\noindent A way to take account of how x- and y-spacings vary jointly is then
to define the two-dimensional spacings as the areas $A_{ij}$ formed by the grid in
Figure \ref{figzone}
$$ \forall i \in \{1,\cdots,n\}, \forall j \in \{1,\cdots,n\}, \quad A_{ij}= D^x_i  D^y_j .$$

\subsection{Uniformity hypothesis}

From now on, I will assume the uniformity hypothesis $H_0$: the point pattern
$U$ is uniformly distributed in $[0,1]^2$.

\noindent Let $(D_1,\cdots,D_{n})$ be the spacings corresponding to a (n-1)
uniform sample on $[0,1]$. Then it is easy to see that

\begin{displaymath}
\left\{ \begin{array}{ll}
\mathcal{L} (D^x_1,\cdots,D^x_n) = \mathcal{L} (D_1,\cdots,D_n), \\
\mathcal{L} (D^y_1,\cdots,D^y_n) = \mathcal{L} (D_1,\cdots,D_n), \\
(D^x_1,\cdots,D^x_n) \perp \! \! \! \perp (D^y_1,\cdots,D^y_n) . \\

\end{array} \right. 
\end{displaymath}

\section{Asymptotic normality of additive functions of spacings}

\subsection{Main result}

Many statistics based on one-dimensional spacings are additive functions
$$ V_n = \sum_{i=1}^n g(n D_i) $$

\noindent for a measurable function $g$.


\noindent For the two-dimensional case consider
$$ V_n^{(2)} = \sum_{i=1}^n \sum_{j=1}^n g(n^2 A_{ij}). $$

\noindent The asymptotic normality of $V_n$ was proved by Beirlant \&
\textit{al}. (1991) using the following distributional equivalence (Moran,
1947)
\begin{eqnarray}
 & (n D_1,\cdots,n D_n) \sim ({\displaystyle
   \frac{E_1}{\bar{E}}},\cdots,{\displaystyle \frac{E_n}{\bar{E}}})
 \label{eq1} \\
 & \nonumber \\
\textrm{where} & (E_1, \cdots, E_n) \textrm{ are independent exponentially
  distributed r.v.'s with mean } 1 \nonumber \\
 \textrm{ and} & \bar{E} = {\displaystyle \frac{1}{n} \sum_{i=1}^n E_i} \, \, . \nonumber
\end{eqnarray}

\noindent I will use the same technique to prove the asymptotic normality of
$V_n^{(2)}$ under $H_0$.

\vspace{0.2cm}

\noindent Introduce $(X_1, \cdots, X_n)$ and $(Y_1, \cdots, Y_n)$ two
samples of independent exponentially distributed r.v.'s with mean $1$.
Denote $\bar{X} = {\displaystyle \frac{1}{n} \sum_{i=1}^n X_i}$  and 

\noindent $\bar{Y}
= {\displaystyle \frac{1}{n} \sum_{j=1}^n Y_j}$.
Then define the statistic
$$ G_n = {\displaystyle \sum_{i=1}^n \sum_{j=1}^n
  g\Big(\frac{X_i}{\bar{X}} \frac{Y_j}{\bar{Y}}\Big)} .$$

\noindent From (\ref{eq1}) , $V_n^{(2)}$ and $G_n$ have same distribution. 

\vspace{0.2cm}

From now on I will assume $g$ satisfies the following, where $\varphi$ and
$\psi$ are measurable functions
\begin{eqnarray} 
g \textrm{ continuous on } \mathbb{R}^+, \qquad \qquad \qquad \qquad \qquad
\qquad  \label{eq2} \\
\mathbb{E}g^2(X_1 Y_1) < \infty, \qquad \qquad \qquad \qquad \qquad \qquad
\label{eq2bis} \\
\forall t_0 \in \mathbb{R}^+, \exists \varphi : \mathbb{R} \rightarrow
\mathbb{R}, \forall n \in \mathbb{N}, \mathbb{E} \varphi (\bar{X} \bar{Y}) <
\infty \qquad \qquad \qquad \nonumber \\
\textrm{ and } \forall t<t_0, \forall x \in \mathbb{R}^+,
|g(tx)|<\varphi(x), \label{eq2ter} \\
\forall t_0 \in \mathbb{R}^+, \exists \psi : \mathbb{R} \rightarrow
\mathbb{R}, \forall n \in \mathbb{N}, \mathbb{E} \psi (\bar{X} \bar{Y}) <
\infty \qquad \qquad \qquad \nonumber \\
\textrm{ and } \forall t>t_0, \forall x \in \mathbb{R}^+,
\Big| \frac{g(tx)}{g(t)} \Big|< \psi (x). \label{eq2tetra}
\end{eqnarray} 

\noindent Denote
\begin{eqnarray*}
\mu & = & \mathbb{E} [g(X_1 Y_1)] ,\\
\eta & = & Cov \big(g(X_1 Y_1),g(X_1 Y_2) \big) = Cov \big(g(X_1 Y_1),g(X_2 Y_1) \big) ,\\
c & = & Cov \big(g(X_1 Y_1),X_1 \big) = Cov \big(g(X_1 Y_1),Y_1 \big) .
\end{eqnarray*}

\noindent To justify my decomposition, I make the following argument.

\noindent Using the same Taylor-expansion as Proschan \& Pyke (1964), it
appears that $n^{-3/2} (G_n-n^2 \mu)$ could be asymptotically equivalent in
distribution to

\begin{eqnarray*} 
\! \! \! \! \! \!   & \! \! \! \! \! \! & \! \! \! \! n^{-3/2} \sum_{i=1}^n \sum_{j=1}^n \{ g(X_i Y_j)  - \mu  - (\bar{X}-1)g'(X_i
Y_j)(X_i Y_j) - (\bar{Y}-1)g'(X_i Y_j)(X_i Y_j) \}  \\
\! \! \! \! \! \!  & \! \! \! \! \! \! = & \! \! \! \! n^{-3/2} \sum_{i=1}^n \sum_{j=1}^n \Big\{ g(X_i Y_j) - \mu  - (X_i-1+Y_j-1) \Big(\frac{1}{n^2}
 \sum_{k=1}^n \sum_{l=1}^n g'(X_k Y_l)(X_k Y_l) \Big) \Big\} .
\end{eqnarray*}

\noindent By partial integration it follows that
$$
 \mathbb{E} g'(X_1 Y_1)(X_1 Y_1) = c \quad \Rightarrow \quad \frac{1}{n^2}
 \sum_{k=1}^n \sum_{l=1}^n g'(X_k Y_l)(X_k Y_l) \xrightarrow[n \rightarrow \infty]{P}
 c .$$

\noindent That is why it seems useful to decompose $G_n$ as follows
\begin{eqnarray*}
{\displaystyle \frac{1}{n^{3/2}}} \big(G_n-n^2 \mu \big) & \! \! \! \! = & \! \! \! \! 
  {\displaystyle \frac{S_n}{n^{3/2}}} + {\displaystyle \frac{R_n}{n^{3/2}}} \\
\textrm{where } S_n & \! \! \! \! = & \! \! \! \! \sum_{i=1}^n \sum_{j=1}^n [g(X_i Y_j) - \mu - c(X_i-1) -c(Y_j-1)] \\
\textrm{and } R_n & \! \! \! \! = & \! \! \! \! \sum_{i=1}^n \sum_{j=1}^n \Big[g\Big(\frac{X_i}{\bar{X}} \frac{Y_j}{\bar{Y}}\Big)
- g(X_i Y_j)\Big] + cn^2(\bar{X}-1) + cn^2(\bar{Y}-1) .
\end{eqnarray*}

\vspace{0.2cm}

\noindent As ${\displaystyle \frac{S_n}{n^{3/2}}}$ is a two-sample U-statistic
with mean $0$ and limiting variance 

\vspace{0.2cm}
\noindent $\sigma^2=2(\eta-c^2)$, one has from Van Der Vaart (1998)
$$ \frac{S_n}{n^{3/2}} \xrightarrow[n \rightarrow \infty]{d}
\mathcal{N} (0,\sigma^2) .$$

\noindent I will prove in the next section that
$$ {\displaystyle \frac{\mathbb{E} (R_n^2)}{n^3}} \xrightarrow[n \rightarrow \infty]{} 0 ,$$

\noindent which will yield the following result.

\begin{theo}
Assume $g$ satisfies (\ref{eq2}), (\ref{eq2bis}), (\ref{eq2ter}) and (\ref{eq2tetra}). Then

$$ \frac{1}{n^{3/2}} \{ G_n - n^2 \mu \} \xrightarrow[n \rightarrow \infty]{d} \mathcal{N} (0,\sigma^2) .$$
\end{theo}

\subsection{Behaviour of the remainder term}

Using the independence of ${\displaystyle \frac{X_i}{\bar{X}}}$ and
$\bar{X}$, as well as the independence of ${\displaystyle
  \frac{Y_j}{\bar{Y}}}$ and $\bar{Y}$, one gets

\begin{eqnarray}
\! \! \! \! \! \! \! \! \! {\displaystyle \frac{\mathbb{E} (R_n^2)}{n^3}} & \!
\! \! \! = & \! \! \! \! T_{1,n} + 2 \, T_{2,n} + T_{3,n} \nonumber
\\
\textrm{where } T_{1,n} & \! \! \! \! = & \! \! \! \! \frac{1}{n} \Big[ \mathbb{E} \Big\{ g^2\Big( {\displaystyle \frac{X_1}{\bar{X}}}
{\displaystyle \frac{Y_1}{\bar{Y}}}\Big) \Big\} -2 \mathbb{E} \Big\{ g\Big( {\displaystyle \frac{X_1}{\bar{X}}}
{\displaystyle \frac{Y_1}{\bar{Y}}}\Big) g(X_1 Y_1) \Big\} + \mathbb{E}
{g^2(X_1 Y_1)} \Big], \label{eq3} \\
 & & \nonumber \\
T_{2,n} & \! \! \! \! \! \! = &  \! \! \! \! \! \! \frac{n-1}{n} \Big[ \mathbb{E} \Big\{ g\Big( {\displaystyle \frac{X_1}{\bar{X}}}
{\displaystyle \frac{Y_1}{\bar{Y}}}\Big) g\Big( {\displaystyle \frac{X_1}{\bar{X}}}
{\displaystyle \frac{Y_2}{\bar{Y}}}\Big) \Big\} -2 \mathbb{E} \Big\{ g\Big( {\displaystyle \frac{X_1}{\bar{X}}}
{\displaystyle \frac{Y_1}{\bar{Y}}}\Big) g(X_1 Y_2) \Big\} \nonumber \\
 & & \qquad \qquad \qquad \qquad \qquad \qquad \qquad \qquad \qquad + \mathbb{E} {g(X_1
  Y_1) g(X_1 Y_2)} \Big],  \label{eq3bis} \\
\textrm{and } T_{3,n} & \! \! \! \! = & \! \! \! \! \frac{(n-1)^2}{n} \Big[ \mathbb{E} \Big\{ g \Big( {\displaystyle
  \frac{X_1}{\bar{X}}} {\displaystyle \frac{Y_1}{\bar{Y}}} \Big) g \Big(
{\displaystyle \frac{X_2}{\bar{X}}} {\displaystyle \frac{Y_2}{\bar{Y}}} \Big)
\Big\} -2 \mathbb{E} \Big\{ g \Big( {\displaystyle \frac{X_1}{\bar{X}}}
{\displaystyle \frac{Y_1}{\bar{Y}}} \Big) g(X_2 Y_2) \Big\} + \mu^2 \Big]
\nonumber \\
 & & \qquad \qquad \qquad \qquad \qquad \qquad \qquad \qquad \qquad \qquad
 \qquad \qquad -2c^2 . \label{eq4}
\end{eqnarray}

\subsubsection{Preliminary results}

\noindent The marginal and bivariate densities of the spacings $D_i$ are given
by Pyke (1965) and lead to

\begin{equation}
\mathbb{E} g^2 \Big( {\displaystyle \frac{X_1}{\bar{X}}}
{\displaystyle \frac{Y_1}{\bar{Y}}} \Big) = {\displaystyle
  \frac{(n-1)^2}{n^2}} \int_0^n \int_0^n g^2(xy) \Big( 1- \frac{x}{n}
\Big)^{n-2} \Big( 1- \frac{y}{n} \Big)^{n-2} dy \, dx, \label{eq5}
\end{equation}

\begin{eqnarray}
\mathbb{E} \Big\{ g \Big( {\displaystyle \frac{X_1}{\bar{X}}}
{\displaystyle \frac{Y_1}{\bar{Y}}} \Big) g \Big( {\displaystyle \frac{X_1}{\bar{X}}}
{\displaystyle \frac{Y_2}{\bar{Y}}} \Big) \Big\} & = & {\displaystyle
  \frac{(n-1)^2 (n-2)}{n^3}} \int_0^n \int_0^n \int_0^{n-y}
g(xy) g(xv) \nonumber \\
 & & \Big( 1- \frac{x}{n}
\Big)^{n-2} \Big( 1- \frac{y+v}{n} \Big)^{n-3} dv \, dy \, dx, \label{eq6bis}
\end{eqnarray} 

\begin{eqnarray}
\mathbb{E} \Big\{ g \Big( {\displaystyle \frac{X_1}{\bar{X}}}
{\displaystyle \frac{Y_1}{\bar{Y}}} \Big) g \Big( {\displaystyle \frac{X_2}{\bar{X}}}
{\displaystyle \frac{Y_2}{\bar{Y}}} \Big) \Big\} & = & {\displaystyle
  \frac{(n-1)^2 (n-2)^2}{n^4}} \int_0^n \int_0^n \int_0^{n-x} \int_0^{n-y}
g(xy) g(uv) \nonumber \\
 & & \Big( 1- \frac{x+u}{n}
\Big)^{n-3} \Big( 1- \frac{y+v}{n} \Big)^{n-3} dv \, du \, dy \, dx . \label{eq6}
\end{eqnarray} 


\noindent Using the independence of ${\displaystyle \frac{X_i}{\bar{X}}}$ and
$\bar{X}$, one finds

\begin{eqnarray*}
\mathbb{E} \Big\{ g \Big( {\displaystyle \frac{X_1}{\bar{X}}}
{\displaystyle \frac{Y_1}{\bar{Y}}} \Big) g(X_1 Y_1) \Big\} & = & \mathbb{E}
\Big[ g \Big( {\displaystyle \frac{X_1}{\bar{X}}}
{\displaystyle \frac{Y_1}{\bar{Y}}} \Big) \mathbb{E} \Big\{
g\Big( \frac{X_1}{\bar{X}} \bar{X} \, \frac{Y_1}{\bar{Y}} \bar{Y} \Big) \Big|
    \frac{X_1}{\bar{X}}, \frac{Y_1}{\bar{Y}} \Big\} \Big] \\
 & = & \frac{(n-1)^2}{n^2} \int_0^n \int_0^n g(xy) \mathbb{E} \{g(x \bar{X} y
 \bar{Y}) \} \nonumber \\
 & & \Big( 1- \frac{x}{n} \Big)^{n-2} \Big( 1- \frac{y}{n} \Big)^{n-2} dy \,
 dx,
\end{eqnarray*}

\begin{eqnarray}
\mathbb{E} \Big\{ g \Big( {\displaystyle \frac{X_1}{\bar{X}}}
{\displaystyle \frac{Y_1}{\bar{Y}}} \Big) g(X_1 Y_2) \Big\} & = & \mathbb{E}
\Big[ g \Big( {\displaystyle \frac{X_1}{\bar{X}}}
{\displaystyle \frac{Y_1}{\bar{Y}}} \Big) \mathbb{E} \Big\{
g\Big( \frac{X_1}{\bar{X}} \bar{X} \, \frac{Y_2}{\bar{Y}} \bar{Y} \Big) \Big|
    \frac{X_1}{\bar{X}}, \frac{Y_1}{\bar{Y}},
    \frac{Y_2}{\bar{Y}} \Big\} \Big] \nonumber \\
 & = & \frac{(n-1)^2 (n-2)}{n^3} \int_0^n \int_0^n \int_0^{n-y} g(xy) \mathbb{E} \{g(x \bar{X} v
 \bar{Y}) \} \nonumber \\
 & & \Big( 1- \frac{x}{n}
\Big)^{n-2} \Big( 1- \frac{y+v}{n} \Big)^{n-3} dv \, dy \, dx, \label{eq6ter}
\end{eqnarray}

\begin{eqnarray}
\! \! \! \mathbb{E} \Big\{ g \Big( {\displaystyle \frac{X_1}{\bar{X}}}
{\displaystyle \frac{Y_1}{\bar{Y}}} \Big) g(X_2 Y_2) \Big\} & \! \! \! \! \! =
& \! \! \! \mathbb{E}
\Big[ g \Big( {\displaystyle \frac{X_1}{\bar{X}}}
{\displaystyle \frac{Y_1}{\bar{Y}}} \Big) \mathbb{E} \Big\{
g\Big( \frac{X_2}{\bar{X}} \bar{X} \, \frac{Y_2}{\bar{Y}} \bar{Y} \Big) \Big|
    \frac{X_1}{\bar{X}}, \frac{X_2}{\bar{X}}, \frac{Y_1}{\bar{Y}},
    \frac{Y_2}{\bar{Y}} \Big\} \Big] \nonumber \\
 & \! \! \! = & \! \! \! \! \! \frac{(n-1)^2 (n-2)^2}{n^4} \int_0^n \int_0^n
 \int_0^{n-x} \! \int_0^{n-y} g(xy) \mathbb{E} \{g(u \bar{X} v \bar{Y}) \}
 \nonumber \\
 & & \! \! \! \Big( 1- \frac{x+u}{n} \Big)^{n-3} \Big( 1- \frac{y+v}{n}
 \Big)^{n-3} dv \, du \, dy \, dx .\label{eq7}
\end{eqnarray}

\noindent The following lemma is also needed.

\begin{lem}
If $g \textrm{ continuous on } \mathbb{R}^{+\star}$ and $\mathbb{E}g^2(X_1 Y_1) < \infty$, then $\forall t \in [0,+\infty[$
$$ \lim_{n \rightarrow \infty} \mathbb{E} g(t \bar{X} \bar{Y}) = g(t). $$
\end{lem}
\noindent \textit{\textbf{Proof}}:

Denote by $f_n(u,t) = {\displaystyle \frac{n^n}{\Gamma (n)}
  \frac{u^{n-1}}{t^n}e^{-nu/t}}$ the common density of $t\bar{X}$ and $t\bar{Y}$.

\begin{eqnarray*}
& & \! \! \! \mathbb{E} g(t \bar{X} \bar{Y}) = \mathbb{E} g(\sqrt{t} \bar{X} \sqrt{t}
\bar{Y}) = \int_0^{\infty} \int_0^{\infty} g(uv) f_n(u,\sqrt{t})
f_n(v,\sqrt{t}) du \, dv \\
 & \! \! \! = & \! \! \! \int \int_D g(uv) f_n(u,\sqrt{t}) f_n(v,\sqrt{t}) du \, dv + \int \int_{\bar{D}}
 g(uv) f_n(u,\sqrt{t}) f_n(v,\sqrt{t}) du \, dv
\end{eqnarray*} 

$$ \textrm{where } D = \{(u,v) \in \mathbb{R}^{+2}; \sqrt{t}/2 < u <
 3\sqrt{t}/2; \sqrt{t}/2 < v < 3\sqrt{t}/2\}. $$

\noindent It is easy to see that:

\noindent $ \sqrt{t} \bar{X} \xrightarrow[n \rightarrow \infty]{P} \sqrt{t} $ and $ \sqrt{t} \bar{Y} \xrightarrow[n \rightarrow \infty]{P} \sqrt{t} $.

\vspace{0.2cm}
\noindent Introduce the function $\varphi: \qquad D \rightarrow \mathbb{R}$

$ \qquad \qquad \qquad \qquad \, \, \, \, \qquad (x,y) \rightarrow \varphi(x,y)=g(xy)$.

\vspace{0.2cm}
\noindent From (\ref{eq2}) the function $\varphi$ is continuous and bounded on
$D$ , so by the Helly-Bray theorem one concludes

$$ \int \int_{D} g(uv) f_n(u,\sqrt{t}) f_n(v,\sqrt{t}) du \, dv = \mathbb{E}
\, \varphi (\sqrt{t} \bar{X}, \sqrt{t} \bar{Y}) \xrightarrow[n \rightarrow
\infty]{} \varphi(\sqrt{t},\sqrt{t}) = g(t) .$$

\vspace{0.2cm}

It remains to prove that $\int \int_{\bar{D}} g(uv) f_n(u,\sqrt{t}) f_n(v,\sqrt{t})
du \, dv \xrightarrow[n \rightarrow \infty]{} 0$.

\noindent $K$ and $\tilde{K}$ are two constants. From Beirlant \& \textit{al} (1991), one has $\forall n \geq 1$
$$ f_n(u,\sqrt{t}) \leq K t^{-1/2} n^{1/2} \Big( \frac{u}{\sqrt{t}}
e^{1-u/\sqrt{t}} \Big)^{n-1} .$$

\noindent So: $(u,v) \in \bar{D} \Rightarrow f_n(u,\sqrt{t}) f_n(v,\sqrt{t})
\xrightarrow[n \rightarrow \infty]{} 0$.

\vspace{0.2cm}

\noindent Moreover: $n \geq n_0 \geq 5 \Rightarrow f_n(u,\sqrt{t})
f_n(v,\sqrt{t}) \leq \tilde{K} f_{n_0}(u,\sqrt{t}) f_{n_0}(v,\sqrt{t})$.

\noindent And taking $m \geq 16$, we get, $\forall u \in [0,\sqrt{t}/2] \cup
[3\sqrt{t}/2,+\infty]$
\begin{eqnarray*}
 & & \Big( \frac{u}{\sqrt{t}} e^{1-u/\sqrt{t}} \Big)^m< e^{-u/\sqrt{t}} \\
 & \Rightarrow & f_{n_0}(u,\sqrt{t}) < K^{1/2} t^{-1/2} n_0^{1/2} e^{-u} \qquad
 \textrm{ if we take } n_0 > 16 \sqrt{t} + 1.
\end{eqnarray*}

\noindent As, from (\ref{eq2bis}), $(u,v) \rightarrow g(uv) e^{-u} e^{-v}  \in
\mathcal{L}^1 (\mathbb{R}^2)$, one gets
$$ \int \int_{\bar{D}} g(uv) f_{n_0}(u,\sqrt{t}) f_{n_0}(v,\sqrt{t}) du \, dv
< \infty . $$

 \noindent Lebesgue's dominated-convergence theorem leads to the conclusion.
 $\blacksquare$

\subsubsection{Behaviour of $T_{1,n}$}

$ \forall (x,y) \in [0,n]^2$,

$ \! \! \! |g^2(xy) (1-x/n)^{n-2}
(1-y/n)^{n-2} | \leq g^2(xy) e^4 e^{-x} e^{-y} \quad \in \mathcal{L}^1
([0,n]^2) \textrm{ from (\ref{eq2bis}).} $

\vspace{0.2cm}

\noindent So applying Lebesgue's dominated-convergence theorem to (\ref{eq5}) leads to
\begin{equation}
\lim_{n \rightarrow \infty} \mathbb{E} g^2 \Big( \frac{X_1}{\bar{X}}
\frac{Y_1}{\bar{Y}} \Big) = \int_0^{\infty} \int_0^{\infty} g^2(xy) e^{-x}
e^{-y} dy \, dx = \mathbb{E} g^2 (X_1 Y_1). \label{eq8}
\end{equation}

By Cauchy-Schwarz inequality one gets
\begin{eqnarray}
\lim_{n \rightarrow \infty} \mathbb{E} \Big\{ \Big| g \Big(
\frac{X_1}{\bar{X}} \frac{Y_1}{\bar{Y}} \Big) g(X_1 Y_1) \Big| \Big\} & \leq &
\{ \mathbb{E} g^2(X_1Y_1) \}^{1/2} \quad \lim_{n \rightarrow \infty} \Big\{
\mathbb{E} g^2 \Big( \frac{X_1}{\bar{X}} \frac{Y_1}{\bar{Y}} \Big)
\Big\}^{1/2} \nonumber \\
 & = & \mathbb{E} g^2 (X_1 Y_1) \quad \textrm{from (\ref{eq8})} \nonumber \\
 & < & \infty . \label{eq8bis}
\end{eqnarray}

\noindent (\ref{eq3}), (\ref{eq8}) and (\ref{eq8bis}) lead to
$$ \lim_{n \rightarrow \infty} T_{1,n} = 0 .$$

\subsubsection{Behaviour of $T_{2,n}$}

$ \forall (x,y) \in [0,n]^2$, $ \forall v \in [0,n-y],$

 $ |g(xy) g(xv) (1-x/n)^{n-2}
(1-(y+v)/n)^{n-3} | \leq |g(xy) g(xv) e^5 e^{-x} e^{-y} e^{-v}| \quad \in \mathcal{L}^1
([0,n]^2 \times [0,n-y]) \textrm{ from (\ref{eq2bis}).} $

\vspace{0.2cm}

\noindent So applying Lebesgue's dominated-convergence theorem to (\ref{eq6bis}) leads to
\begin{equation}
\lim_{n \rightarrow \infty} \mathbb{E} g \Big( \frac{X_1}{\bar{X}}
\frac{Y_1}{\bar{Y}} \Big) g \Big( \frac{X_1}{\bar{X}}
\frac{Y_2}{\bar{Y}} \Big) = \mathbb{E} g(X_1 Y_1) g(X_1 Y_2). \label{eq8ter}
\end{equation}

\noindent Introduce the function $h_n: \mathbb{R}^{+*} \rightarrow \mathbb{R}$

$ \qquad \quad \qquad \qquad \, \, \, \, \, \qquad \qquad t \rightarrow h_n(t)=\mathbb{E} g(t \bar{X} \bar{Y})$.

\noindent Lemma $1$ gives
$$g(xy) h_n(xv) (1-x/n)^{n-2} (1-(y+v)/n)^{n-3} \xrightarrow[n \rightarrow
\infty]{} g(xy) g(xv) e^{-x} e^{-y} e^{-v} .$$

Denote $t_1 \in \mathbb{R}^{+*}$. From (\ref{eq2tetra}), one gets

$ \exists \psi : \mathbb{R} \rightarrow
\mathbb{R}, \forall t>t_1, \forall x \in \mathbb{R}^+, \Big| {\displaystyle
\frac{g(txy)}{g(t)} \Big|} < \psi(x)$

\begin{eqnarray*}
\Rightarrow \Big| \frac{h_n(t)}{g(t)} \Big| & = & \Big| \int_0^{\infty} \int_0^{\infty}
\frac{g(txy)}{g(t)} f_n(x,1) f_n(y,1) dy \, dx \Big| \\
 & < & \int_0^{\infty} \int_0^{\infty}
 \psi(xy) f_n(x,1) f_n(y,1) dy  \, dx = \mathbb{E} \psi (\bar{X} \bar{Y}) = \Psi(n).
\end{eqnarray*}

So, $\forall (x,v) \in \mathbb{R}^{+2}, xv>t_1 \Rightarrow |g(xy) h_n(xv)
(1-x/n)^{n-2} (1-(y+v)/n)^{n-3}| < |g(xy) g(xv) \Psi(n) e^5 e^{-x} e^{-y} e^{-v}| \in \mathcal{L}^1([0,n]^2 \times [0,n-y]).$

\vspace{0.2cm}

From (\ref{eq2ter}), one gets

$ \exists \varphi : \mathbb{R} \rightarrow
\mathbb{R}, \forall t<t_1, \forall x \in \mathbb{R}^+, |g(tx)|<\varphi(x)$

\noindent $\Rightarrow \forall t<t_1, |h_n(t)| < \mathbb{E} \varphi (\bar{X}
\bar{Y}) = \Phi(n).$

So, $\forall (x,v) \in \mathbb{R}^{+2}, xv<t_1 \Rightarrow |g(xy) h_n(xv) (1-x/n)^{n-2} (1-(y+v)/n)^{n-3}| < |g(xy)
\Phi(n) e^5 e^{-x} e^{-y} e^{-v}| \in \mathcal{L}^1([0,n]^2 \times [0,n-y]).$ 

\vspace{0.2cm}

\noindent So applying Lebesgue's dominated-convergence theorem to (\ref{eq6ter}) leads to
\begin{equation}
\lim_{n \rightarrow \infty} \mathbb{E} \Big\{  g \Big(
\frac{X_1}{\bar{X}} \frac{Y_1}{\bar{Y}} \Big) g(X_1 Y_2) \Big\} = \mathbb{E}
g(X_1 Y_1) g(X_1 Y_2). \label{eq9}
\end{equation}

\vspace{0.2cm}

\noindent (\ref{eq3bis}), (\ref{eq8ter}) and (\ref{eq9}) lead to
$$ \lim_{n \rightarrow \infty} T_{2,n} = 0 .$$

\subsubsection{Behaviour of $T_{3,n}$}

\noindent From (\ref{eq6}), one gets
\begin{eqnarray}
& & \! \! \! \! \! \! \! \! \! \! \! \frac{(n-1)^2}{n} \mathbb{E} \Big\{ g \Big( {\displaystyle \frac{X_1}{\bar{X}}}
{\displaystyle \frac{Y_1}{\bar{Y}}} \Big) g \Big( {\displaystyle
  \frac{X_2}{\bar{X}}} {\displaystyle \frac{Y_2}{\bar{Y}}} \Big) \Big\}
\nonumber \\
& = & 
{\displaystyle \frac{(n-1)^4 (n-2)^2}{n^5}} \int_0^{\infty} \int_0^{\infty}
\int_0^{\infty} \int_0^{\infty} g(xy) g(uv) \Big( 1+ {\displaystyle
  \frac{3(x+u) -(x+u)^2/2}{n}}\Big) \nonumber \\
& & \Big( 1+ {\displaystyle \frac{3(y+v)
 -(y+v)^2/2}{n}}\Big) e^{-x} e^{-y} e^{-u} e^{-v}  dv \, du \, dy \, dx +I_n
\nonumber \\
 & = & {\displaystyle \frac{(n-1)^4 (n-2)^2}{n^5}} \mu^2 +12 \mu {\displaystyle
  \frac{(n-1)^4 (n-2)^2}{n^6}} \mathbb{E} [X_1 g(X_1 Y_1)] +
\mathcal{O}(n^{-1}) \label{eq10} \\
 & - & 2 \mu {\displaystyle \frac{(n-1)^4 (n-2)^2}{n^6}} \mathbb{E} [X_1^2
 g(X_1 Y_1)] - 2 {\displaystyle \frac{(n-1)^4 (n-2)^2}{n^6}} \mathbb{E}^2 [X_1
 g(X_1Y_1)] + I_n \nonumber 
\end{eqnarray}

\noindent where $I_n=$
\begin{eqnarray*}
 & & \! \! \! \! \! \! \! \! \! \! \! \! {\displaystyle \frac{(n-1)^4 (n-2)^2}{n^5}} \Big\{ \int_0^n
  \int_0^n \int_0^{n-x} \int_0^{n-y} g(xy) g(uv) \Big( 1-{\displaystyle
    \frac{x+u}{n}} \Big)^{n-3} \Big( 1-{\displaystyle \frac{y+v}{n}}
  \Big)^{n-3} \\
 & & \! \! \! \! \! \! \! \! \! \! \! \! dv \, du \, dy \, dx -
 \int_0^{\infty} \int_0^{\infty} \int_0^{\infty} \int_0^{\infty} g(xy) g(uv)
 \Big( 1+ {\displaystyle \frac{3(x+u) -(x+u)^2/2}{n}}\Big) \\
 & & \! \! \! \! \! \! \! \! \! \! \! \! \Big(
  1+ {\displaystyle \frac{3(y+v) -(y+v)^2/2}{n}}\Big) e^{-x} e^{-y} e^{-u}
  e^{-v} dv \, du \, dy \, dx \Big\} .
\end{eqnarray*}

\noindent From (\ref{eq7}), one gets
\begin{eqnarray*}
U_n & \! \! \! = & \! \! \! \frac{(n-1)^2}{n} \mathbb{E} \Big\{ g \Big( {\displaystyle \frac{X_1}{\bar{X}}}
{\displaystyle \frac{Y_1}{\bar{Y}}} \Big) g (X_2 Y_2) \Big\} \\
& \! \! \! = & \! \! \! \frac{(n-1)^4 (n-2)^2}{n^5} \int_0^n \int_0^n \int_0^{n-x} \int_0^{n-y}
 \Big( \int_0^{\infty} \int_0^{\infty} g(ab) f_n(a,u) f_n(b,v) da \, db
 \Big) \\
& \! \! \! & \! \! \! g(xy) \Big( 1-{\displaystyle
    \frac{x+u}{n}} \Big)^{n-3} \Big( 1-{\displaystyle \frac{y+v}{n}}
  \Big)^{n-3} dv \, du \, dy \, dx \\
& \! \! \! = & \! \! \! \frac{(n-1)^4 (n-2)^2}{n^5} \int_0^n \int_0^n \int_0^{\infty}
\int_0^{\infty}
 g(xy) \Big\{ \int_0^{n-x} f_n(a,u) \Big( 1-{\displaystyle
    \frac{x+u}{n}} \Big)^{n-3} du \Big\} \\
& \! \! \! & \! \! \! g(ab) \Big\{ \int_0^{n-y} f_n(b,v) \Big(
  1-{\displaystyle \frac{y+v}{n}} \Big)^{n-3} dv \Big\} db \, da \, dy \, dx .
\end{eqnarray*}

\noindent Now the first integral in braces is equal to 
$$ \frac{n}{n-1} e^{-a} \Big( 1-\frac{x}{n} \Big)^{n-3} e^{-ax/(n-x)} \bigg( 1+
\frac{na}{(n-x)(n-2)} \bigg) $$

\begin{eqnarray*}
\! \! \! \Rightarrow U_n & \! \! \! = & \! \! \! \frac{(n-1)^2 (n-2)^2}{n^3} \int_0^n \int_0^n \int_0^{\infty}
 \int_0^{\infty} g(xy) g(ab) e^{-a} \Big( 1-\frac{x}{n} \Big)^{n-3} \\
 & & \! \! \! e^{-ax/(n-x)} \bigg( 1+ \frac{na}{(n-x)(n-2)} \bigg) e^{-b} \Big(
 1-\frac{y}{n} \Big)^{n-3} e^{-ay/(n-y)} \\
 & & \! \! \! \bigg( 1+ \frac{nb}{(n-y)(n-2)} \bigg) db \, da \, dy
 \, dx \\
 & \! \! \! = & \! \! \! \frac{(n-1)^2 (n-2)^2}{n^3} \int_0^n \int_0^n \int_0^{\infty}
 \int_0^{\infty} g(xy) g(ab) e^{-a} \Big( 1-\frac{x}{n} \Big)^{n-3} \\
 & & e^{-b}
 \Big( 1-\frac{y}{n} \Big)^{n-3} \bigg( 1-\frac{ax}{n-x} +
 \frac{na}{(n-x)(n-2)} \bigg) \\
 & & \! \! \! \bigg( 1-\frac{by}{n-y} +
 \frac{nb}{(n-y)(n-2)} \bigg) db \, da \, dy \, dx + J_n
\end{eqnarray*}

\noindent where $J_n =$
\begin{eqnarray*}
 & & \! \! \! \frac{(n-1)^2 (n-2)^2}{n^3} \int_0^n \int_0^n \int_0^{\infty}
 \int_0^{\infty} g(xy) g(ab) e^{-a} \Big( 1-\frac{x}{n} \Big)^{n-3} e^{-b}
 \Big( 1-\frac{y}{n} \Big)^{n-3} \\
 & & \! \! \! \bigg\{ e^{-ax/(n-x)} e^{-by/(n-y)} \bigg( 1+ \frac{na}{(n-x)(n-2)} \bigg)
 \bigg( 1+ \frac{nb}{(n-y)(n-2)} \bigg) \\
 & & \! \! \! - \bigg( 1-\frac{ax}{n-x} +
 \frac{na}{(n-x)(n-2)} \bigg) \bigg( 1-\frac{by}{n-y} +
 \frac{nb}{(n-y)(n-2)} \bigg) \bigg\} db \, da \, dy \, dx . 
\end{eqnarray*}

\noindent Substituting ${\displaystyle \frac{ax}{n} \Big( 1- \frac{x}{n}
  \Big)^{-1}}$ to ${\displaystyle \frac{ax}{n-x}}$, one gets that $U_n$ equals

\begin{eqnarray}
 & & \! \! \! \frac{(n-1)^2 (n-2)^2}{n^3} \mu^2 + \frac{6(n-1)^2 (n-2)^2}{n^4} \mu
 \mathbb{E} [X_1 g(X_1 Y_1)]  \nonumber \\
 & \! \! \! - & \! \! \! \frac{2(n-1)^2 (n-2)^2}{n^4} \mathbb{E} [X_1 g(X_1 Y_1)] \int_0^n \int_0^n
 x g(xy) \Big(1-\frac{x}{n}\Big)^{n-4} \Big(1-\frac{y}{n}\Big)^{n-3} dy \, dx
 \nonumber \\
 & \! \! \! + & \! \! \! \frac{2(n-1)^2 (n-2)}{n^3} \mathbb{E} [X_1 g(X_1 Y_1)] \int_0^n \int_0^n
 g(xy) \Big(1-\frac{x}{n}\Big)^{n-4} \Big(1-\frac{y}{n}\Big)^{n-3} dy \, dx \nonumber \\
 & \! \! \! - & \! \! \! \frac{(n-1)^2 (n-2)^2}{n^4} \mu \mathbb{E} [X_1^2
 g(X_1 Y_1)] + J_n  + K_n + \mathcal{O}(n^{-1}) \label{eq11}
\end{eqnarray}

\noindent where $K_n=$
\begin{eqnarray*}
 & & \frac{(n-1)^2 (n-2)^2}{n^3} \mu \Big[ \int_0^n \int_0^n g(xy) \Big(
1-\frac{x}{n} \Big)^{n-3} \Big( 1-\frac{y}{n} \Big)^{n-3} dy \, dx \\
 & & - \int_0^{\infty} \int_0^{\infty} g(xy) e^{-x} e^{-y} \Big( 1+ \frac{3x
   -x^2/2}{n} \Big) \Big( 1+ \frac{3y-y^2/2}{n} \Big) dy \, dx \Big] .
\end{eqnarray*}

\noindent So, from (\ref{eq4}), (\ref{eq10}) et (\ref{eq11}) we get
\begin{equation}
 T_{3,n} = A_n + \mathcal{O}(n^{-1}) + I_n + J_n + K_n - 2c^2 \label{eq12} 
\end{equation}

\vspace{-0.8cm}
\begin{eqnarray*}
\textrm{where } A_n & \! \! \! = & \! \! \! (n-8) \mu^2 +12\mu \mathbb{E} [X_1 g(X_1
Y_1)] -2\mu \mathbb{E} [X_1^2 g(X_1 Y_1)] -2\mathbb{E}^2 [X_1 g(X_1 Y_1)] \\
 & \! \! \! - & \! \! \! 2(n-6) \mu^2 -12\mu \mathbb{E} [X_1
 g(X_1 Y_1)]  +2\mu \mathbb{E} [X_1^2 g(X_1 Y_1)]  + (n-2) \mu^2 \\
 & \! \! \! + & \! \! \! 4\mathbb{E} [X_1 g(X_1 Y_1)] \int_0^n \int_0^n
 x g(xy) \Big(1-\frac{x}{n}\Big)^{n-4} \Big(1-\frac{y}{n}\Big)^{n-3} dy \, dx \\
 & \! \! \! - & \! \! \! 4 \mathbb{E} [X_1 g(X_1 Y_1)] \int_0^n \int_0^n
 g(xy) \Big(1-\frac{x}{n}\Big)^{n-4} \Big(1-\frac{y}{n}\Big)^{n-3} dy \, dx + \mathcal{O}(n^{-1})
\end{eqnarray*}

\vspace{-0.6cm}
\begin{eqnarray}
 \Rightarrow A_n & \xrightarrow[n \rightarrow \infty]{} & 2\mu^2 -4\mu \mathbb{E} [X_1 g(X_1 Y_1)] +2 \mathbb{E}^2 [X_1 g(X_1
 Y_1)] \nonumber \\
 \Rightarrow A_n & \xrightarrow[n \rightarrow \infty]{} & 2 c^2 . \label{eq13}
\end{eqnarray}

\vspace{0.2cm}
It now suffices to show that $I_n + J_n + K_n = o(1)$.

\noindent A Taylor-expansion leads to

\vspace{-0.4cm}
\begin{equation} 
\forall x \in [1,\sqrt{n}]:
 \Big( 1- \frac{x}{n} \Big)^{n-3} = e^{-x} \Big[ 1 + \frac{3x
  - x^2/2}{n} + \mathcal{O}\big(n^{-2} (x^2+x^4)\big) \Big].  \label{eq14} 
\end{equation}

\noindent Denote $b_n=4 \log{n}$. Choosing $n$ large enough gives
\begin{equation}
\forall x \in [b_n,n], \Big( 1- {\displaystyle \frac{x}{n}} \Big)^{n-3} <
e^{-x}. \label{eq15}
\end{equation}

\noindent By Cauchy-Schwarz inequality one gets
\vspace{-0.2cm}
\begin{eqnarray*}
\frac{n^3}{(n-1)^2 (n-2)^2 \mu} |K_n| & \! \! \! = & \! \! \! \Big| \int_0^{\infty}
\int_0^{\infty} g(xy) \Big[ \Big\{ \Big(
1-\frac{x}{n} \Big)^+ \Big( 1-\frac{y}{n} \Big)^+ \Big\}^{n-3} \\
 & \! \! \! - & \! \! \! e^{-x} e^{-y} \Big( 1+ \frac{3x
   -x^2/2}{n} \Big) \Big( 1+ \frac{3y-y^2/2}{n} \Big) \Big] dy \, dx \Big| \\
 & \! \! \! \leq & \! \! \! [\mathbb{E}g^2(X_1 Y_1)]^{1/2} \Bigg[ \int_0^{\infty} \int_0^{\infty}
 \Big[ \Big\{ \Big(
1-\frac{x}{n} \Big)^+ \Big( 1-\frac{y}{n} \Big)^+ \Big\}^{n-3} \\
 & \! \! \! - & \! \! \! e^{-x} e^{-y} \Big( 1+ \frac{3x
   -x^2/2}{n} \Big) \Big( 1+ \frac{3y-y^2/2}{n} \Big) \Big]^2 e^x e^y dy \, dx \Bigg]^{1/2}
\end{eqnarray*}
\vspace{-0.5cm}
where $x^+ =x$ if $x>0$, $0$ elsewhere. 

\vspace{0.5cm}

\noindent Denote $E_n$ as the double integral 
\begin{eqnarray*}
 & & \! \! \! \! \! \int_0^{\infty} \int_0^{\infty} \Big[ \Big\{ \Big( 1-\frac{x}{n} \Big)^+
\Big( 1-\frac{y}{n} \Big)^+ \Big\}^{n-3} - e^{-x} e^{-y} \\
 & & \Big( 1+ \frac{3x -x^2/2}{n} \Big) \Big( 1+ \frac{3y-y^2/2}{n} \Big)
 \Big]^2 e^x e^y dy \, dx.
\end{eqnarray*}

\noindent Using (\ref{eq14}) and (\ref{eq15}) and following the technique of
Does \& Klaassen (1984), one can prove: $E_n = \mathcal{O}(n^{-4})$ 
\vspace{-0.2cm}
\begin{equation}
\Rightarrow K_n = \mathcal{O}(n^{-1}). \label{eq16}
\end{equation}

\noindent The same arguments are used for $I_n$
\vspace{-0.2cm}
\begin{equation}
\Rightarrow I_n = \mathcal{O}(n^{-1}) .\label{eq17}
\end{equation}

\noindent Using the inequalities $1-z \geq e^{-z} \geq 1-z+z^2/2$, $z \geq 0$,
the expression in braces in the definition of $J_n$ is bounded below by

\vspace{-0.4cm}
\begin{eqnarray*}
 & & - \frac{a^2nx}{(n-x)^2(n-2)} - \frac{n^2a^2bx}{(n-2)^2
   (n-x)^2(n-y)} \\ 
 & - & \frac{b^2ny}{(n-y)^2(n-2)} - \frac{n^2b^2ay}{(n-2)^2
   (n-y)^2(n-x)}
\end{eqnarray*}

\noindent and bounded above by

\vspace{-0.4cm}
\begin{eqnarray*}
  & & \frac{a^2x^2}{2(n-x)^2} +
 \frac{na^2bx^2}{2(n-2)(n-x)^2(n-y)} \\
 & + & \frac{a^2b^2x^2y^2}{4(n-x)^2(n-y)^2} +
 \frac{na^2b^3x^2y^2}{4(n-2)(n-x)^2(n-y)^3} \\
 & + & \frac{na^3x^2}{2(n-2)(n-x)^3} + \frac{n^2a^3bx^2}{2(n-2)^2(n-x)^3(n-y)}
 \\
 & + & \frac{na^3b^2x^2y^2}{4(n-2)(n-x)^3(n-y)^2} +
 \frac{n^2a^3b^3x^2y^2}{4(n-2)^2(n-x)^3(n-y)^3} \\
 & + & \frac{b^2y^2}{2(n-y)^2} +
 \frac{nb^2ay^2}{2(n-2)(n-y)^2(n-x)} \\
 & + & \frac{b^2a^2y^2x^2}{4(n-y)^2(n-x)^2} +
 \frac{nb^2a^3y^2x^2}{4(n-2)(n-y)^2(n-x)^3} \\
 & + & \frac{nb^3y^2}{2(n-2)(n-y)^3} + \frac{n^2b^3ay^2}{2(n-2)^2(n-y)^3(n-x)}
 \\
 & + & \frac{nb^3a^2y^2x^2}{4(n-2)(n-y)^3(n-x)^2} +
 \frac{n^2b^3a^3y^2x^2}{4(n-2)^2(n-y)^3(n-x)^3} \, .
\end{eqnarray*}

\noindent Hence one gets

\vspace{-0.4cm}

\begin{eqnarray*} 
 \frac{n^3}{(n-1)^2 (n-2)^2} |J_n| & \leq & \frac{1}{2n^2} \int_0^n
\int_0^n \int_0^{\infty} \int_0^{\infty} |g(xy) g(ab)| a^2 x^2 e^{-a} e^{-b} \\
 & & \Big( 1-\frac{x}{n} \Big)^{n-5} \Big( 1-\frac{y}{n} \Big)^{n-3} db \, da \, dy\, dx \\
 & & \\
 & + & \cdots \\
 & + & \frac{1}{n(n-2)^2} \int_0^n
\int_0^n \int_0^{\infty} \int_0^{\infty} |g(xy) g(ab)| a b^2 y e^{-a} e^{-b}\\
 & & \Big( 1-\frac{x}{n} \Big)^{n-4} \Big(
 1-\frac{y}{n} \Big)^{n-5} db \, da \, dy \, dx
\end{eqnarray*}

\begin{equation}
\Rightarrow J_n = \mathcal{O}(n^{-1}) . \label{eq18}
\end{equation}

\vspace{0.2cm}

\noindent (\ref{eq12}), (\ref{eq13}), (\ref{eq16}), (\ref{eq17}) and
(\ref{eq18}) lead to
$$ \lim_{n \rightarrow \infty} T_{3,n} = 0 .$$

\section*{References}

\begin{description}

\item {\sc Beirlant,} J., {\sc Janssen,} P. and {\sc Veraverbeke,} N. (1991).
\newblock On the asymptotic normality of functions of uniform spacings. 
\newblock {\it The Canadian Journal of Statistics}, \textbf{19}, 93-101.

\item {\sc Cressie,} N.A.C. (1993).
\newblock {\it Statistics for Spatial Data}.
\newblock Wiley, New York.

\item {\sc Cucala,} L. and {\sc Thomas-Agnan,} C. (2005).
\newblock Tests for spatial randomness based on spacings. 
\newblock {\it Submitted to The Canadian Journal of Statistics}.

\item {\sc Deheuvels,} P. (1983a).
\newblock Spacings and applications.
\newblock {\it Probability and Statistical Decision Theory. Volume A} (F.
Konecny, J. Mogyor\'odi and W. Wertz, \textit{editors}). Reidel, Dordrecht, 1-30.

\item {\sc Deheuvels,} P. (1983b).
\newblock Strong bounds for multidimensional spacings.
\newblock {\it Zeitschrift f\"ur Wahrscheinlichkeitstheorie und Verwandte Gebiete}, \textbf{64}, 411-424.

\item {\sc Does,} R.J.M.M. and {\sc Klaassen,} C.A.J. (1984).
\newblock The Berry-Esseen theorem for functions of uniform spacings.
\newblock {\it Zeitschrift f\"ur Wahrscheinlichkeitstheorie und Verwandte Gebiete}, \textbf{65}, 461-472.

\item {\sc Janson,} S. (1987).
\newblock Maximal spacings in several dimensions.
\newblock {\it Annals of Probability}, \textbf{15}, 274-280.

\item {\sc Justel,} A., {\sc Pe\~na,} D. and {\sc Zamar,} R. (1997).
\newblock A multivariate Kolmogorov-Smirnov test of goodness of fit.
\newblock {\it Statistics and Probability Letters}, \textbf{35}, 251-259.

\item {\sc Le Cam,} L. (1958).
\newblock Un théorème sur la division d'un intervalle par des points pris au
hasard.
\newblock {\it Publications de l'Institut de Statistique de l'Université de Paris}, \textbf{7}, 7-16.

\item {\sc Moller,} J. and {\sc Waagepetersen,} R.P. (2004).
\newblock {\it Statistical Inference and Simulation for Spatial Point
  Processes}.
\newblock Chapman \& Hall/CRC, Boca Raton, Florida.

\item {\sc Moran,} P.A.P. (1947).
\newblock The random division of an interval.
\newblock {\it Journal of the Royal Statistical Society Series B}, \textbf{9}, 92-98.

\item {\sc Proschan,} F. and {\sc Pyke,} R. (1964).
\newblock Asymptotic normality of certain test statistics of exponentiality.
\newblock {\it Biometrika}, \textbf{51}, 253-256.

\item {\sc Pyke,} R. (1965).
\newblock Spacings.
\newblock {\it Journal of the Royal Statistical Society Series B}, \textbf{27}, 395-449.


\item {\sc Van Der Vaart,} A.W. (1998).
\newblock {\it Asymptotic Statistics}.
\newblock Cambridge University Press.

\item {\sc Zimmerman,} D.L. (1993).
\newblock A Bivariate Cramer-Von Mises Type of Test for Spatial Randomness.
\newblock {\it Applied Statistics}, \textbf{42}, 43-54.

\end{description}

\end{document}